\documentclass[10pt]{article}
\usepackage{amsmath}
\usepackage{amssymb}
\usepackage{mathrsfs}
\usepackage{amsfonts}
\usepackage{mathrsfs,amscd,amssymb,amsthm,amsmath,bm,graphicx}
\usepackage{graphicx}
\usepackage{cite}
\usepackage{color}
\usepackage[shortcuts]{extdash}
\makeatletter

\renewcommand{\@seccntformat}[1]{{\csname the#1\endcsname}{\normalsize.}\hspace{.5em}}
\makeatother

\def \[{\begin{equation}}
\def \]{\end{equation}}

\newtheorem{thm}{Theorem}[section]

\newtheorem{lem}[thm]{Lemma}

\begin{document}
\setlength{\baselineskip}{13pt}
\begin{center}{\Large \bf The Laplacians, Kirchhoff index and complexity of linear M\"{o}bius and cylinder octagonal-quadrilateral networks
}

\vspace{4mm}

{\large Jia-Bao Liu $^{1,*}$, Lu-Lu Fang $^{1,*}$, Qian Zheng $^{1}$, Xin-Bei Peng $^{1}$}\vspace{2mm}

{\small $^{1}$School of Mathematics and Physics, Anhui Jianzhu University, Hefei 230601, P.R. China\\}
\vspace{2mm}
\end{center}

\footnotetext{E-mail address: liujiabaoad@163.com, fangluluajd@163.com,zhengqian19960202@163.com, pengxinbeiajd@163.com.}

\footnotetext{* Corresponding author.}

{\noindent{\bf Abstract.}\ \ Spectrum graph theory not only facilitate comprehensively reflect the topological structure and dynamic characteristics of networks, but also offer significant and noteworthy applications in theoretical chemistry, network science and other fields. Let $L_{n}^{8,4}$ represent a linear octagonal-quadrilateral network, consisting of $n$ eight-member ring and $n$ four-member ring. The M\"{o}bius graph $Q_{n}(8,4)$ is constructed by reverse identifying the opposite edges, whereas cylinder graph $Q'_{n}(8,4)$ identifies the opposite edges by order. In this paper, the explicit formulas of Kirchhoff indices and complexity of $Q_{n}(8,4)$ and $Q'_{n}(8,4)$ are demonstrated by Laplacian characteristic polynomials according to decomposition theorem and Vieta's theorem. In surprise, the Kirchhoff index of $Q_{n}(8,4)$\big($Q'_{n}(8,4)$\big) is approximately one-third half of its Wiener index as $n\to\infty$.\vspace{2mm}
\section{Introduction}\label{sct1}
\ \ \ \ \ The graphs in this paper are simple, undirected and connected. Firstly recall some definitions that most used in graph theory. Suppose $G$ as a simple undirected graph with $\big|V_G\big|=n$ and $\big|E_G\big|=m$. The adjacent matrix $A(G)=[a_{i,j}]_{n\times n}$ of $G$ is a symmetric matrix. When the vertex $i$ is adjacent to $j$, we define $a_{i,j}=1$ and $a_{i,j}=0$ otherwise. For a vertex $i\in V_G$, $d_{i}$ is the number of edges arised from $i$ and the diagonal matrix $D(G)=diag\{d_{1},d_{2},\cdots,d_{n}\}$ is degree matrix. Subsequently, the Laplacian matrix is specified as $L(G)=D(G)-A(G)$, which spectrum can be expressed by $0=\mu_{1}<\mu_{2}\leq\cdots\leq\mu_{n}$. For more notations, one can be referred to \cite{F.R}.

Alternatively, The $(m,n)$-entry of the Laplacian matrix can be noted by
\begin{eqnarray}
\big(L(G)\big)_{mn}=
\begin{cases}
d_{m}, & m=n;\\
-1, & m\neq{n}~and~ v_{m}\backsim v_{n};\\
0, & otherwise.
\end{cases}
\end{eqnarray}

The traditional distance for every pairs of vertices $i$ and $j$ is $d_{G}(v_i,v_j)=\{d_{ij}|\ the\ shortset\ length\ of\ v_{i}\ \\ and\ v_{j}\}$. For many parameters used to describe the structure of a graph in chemical and mathematical researches\cite{B,C.,D}, the most commonly used one is the distance-based parameter Wiener index\cite{Wiener,A.D}, which is known to
\begin{eqnarray*}
W(G)=\sum_{i<j}d_{ij}.
\end{eqnarray*}

Suppose that each edge of a connected graph $G$ is regarded as a unit resistor, and the resistance distance\cite{D.J} between any two vertices $i$ and $j$ is recorded as $r_{ij}$. Similar to Wiener index, according to the resistance distance, the expression of Kirchhoff index \cite{D.,D.J.} are given, namely
\begin{eqnarray*}
Kf(G)=\sum_{i<j}r_{ij}.
\end{eqnarray*}

It's well known that scholars in various fields have a strong interest in the study of the Kirchhoff index, this has promoted researchers to try some calculation methods to compute the Kirchhoff index of a given graph and get its closed formula. As early as 1985, Y.L. Yang et al.\cite{Y.L} gave the corresponding decomposition theorem for the corresponding matrix of linear structure. Up to now, many scholars have used this method to describe the Kirchhoff index of a series of linear hydrocarbons chain. For instances, in 2007, Yang et al.\cite{Y.} calculated the Kirchhoff index of linear hexagonal chain by using the decomposition of Laplacian matrix. In 2019, Geng et al.\cite{G.W} used this method to calculate the M\"{o}bius chain and cylinder chain of phenylenes chain. Shi and Liu et al.\cite{Z.L} computed the Kirchhoff index of linear octagonal-quadrilateral network in 2020. For more information, see \cite{H.B,E.,F.,G.,H.,M.,K.}. After learning the excellent works of scholars, this paper computes the Kirchhoff indices, Wiener indices and the complexity of M\"{o}bius graph and its cylinder graph of linear octagonal-quadrilateral network.\\

Let $L_{n}^{8,4}$ be the linear octagonal-quadrilateral network as illustrated in Figure 1, and octagons and quadrilaterals are connected by a common edge. Then the corresponding M\"{o}bius graph $Q_{3}\big(8,4\big)$ of octagonal-quadrilateral network is obtained by the reverse identification of the opposite edge by $L_{n}^{8,4}$, and its cylinder graph $Q'_{3}\big(8,4\big)$ of octagonal-quadrilateral network is obtained by identifying the opposite edge of $L_{n}^{8,4}$. An illustration of this are given in Figure 2. Obviously, we can obtained that $\big|V_{Q_{n}}(8,4)\big|=8n$,~$\big|E_{Q_{n}}(8,4)\big|=10n$ and $\big|V_{Q'_{n}}(8,4)\big|=8n$,~ $\big|E_{Q'_{n}}(8,4)\big|=10n.$
\begin{figure}[htbp]
\centering\includegraphics[width=16.5cm,height=4cm]{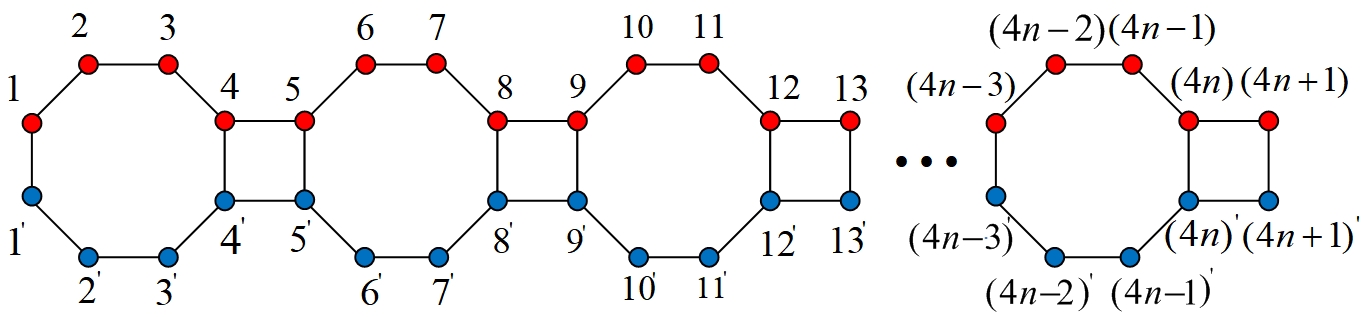}
\caption{Linear octagonal-quadrilateral networks. }
\end{figure}
\begin{figure}[htbp]
\centering\includegraphics[width=13cm,height=4cm]{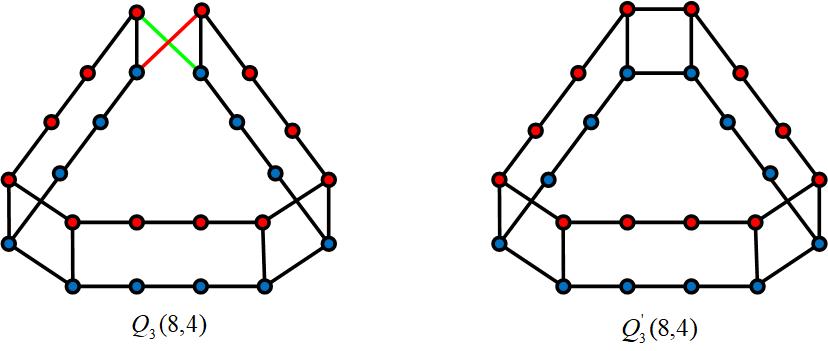}
\caption{A special class of octagonal-quadrilateral networks for $n=3$.}
\end{figure}

The rest of this paper is organized as follows: In Section 2, we introduce some basic notations and related lemmas. In Section 3, the Kirchhoff index, the Wiener index and the complexity of $Q_{n}\big(8,4\big)$ and $Q'_{n}\big(8,4\big)$ are be computed. In Section 4, we demonstrate the ratio of Wiener index and Kirchhoff index of $Q_{n}\big(8,4\big)$ and $Q'_{n}\big(8,4\big)$. Finally, we provide insights into the subsequent extensions in Section 5.
\section{Preliminary}\label{sct2}
\ \ \ \ \ We recommend the notion of symbols and related calculation methods. This in turn would help to the proof of the rest of the article. The characteristic polynomial of matrix $A$ of order $n$ is defined as $P_{A}(x)=det\big(xI-A\big)$. We note that $\varpi$ is an automorphism of $G$, we can write the product of disjoint 1-cycles and transposition, namely
\begin{equation*}
\mathscr{\varpi}=(\bar{1})(\bar{2})\cdots(\bar{m})(1,1')(2,2')\cdots(k,k').
\end{equation*}

Then $\big|V(G)\big|=m+2k$,~ let $v_{0}=\big\{\bar{1},\bar{2},\cdots \bar{m}\big\},~ v_{1}=\big\{1,2\cdots k\big\},~v_{2}=\big\{1',2'\cdots k'\big\}.$ Thus the Laplacian matrix can be formed by block matrix, then
\begin{equation}
L(G)=\left(
\begin{array}{ccc}
L_{V_{0}V_{0}}& L_{V_{0}V_{1}}& L_{V_{0}V_{2}}\\
L_{V_{1}V_{0}}& L_{V_{1}V_{1}}& L_{V_{1}V_{2}}\\
L_{V_{2}V_{0}}& L_{V_{2}V_{1}}& L_{V_{2}V_{2}}\\
\end{array}
\right),
\end{equation}
where
\begin{equation}
L_{V_{0}V_{1}}=L_{V_{0}V_{2}},~ L_{V_{1}V_{2}}=L_{V_{2}V_{1}},~and~ L_{V_{1}V_{1}}=L_{V_{2}V_{2}}.
\end{equation}

Let
\begin{equation}
P=\left(
  \begin{array}{ccc}
  I_{m}& 0& 0\\
  0& \frac{1}{\sqrt{2}}I_{k}& \frac{1}{\sqrt{2}}I_{k}\\
  0& \frac{1}{\sqrt{2}}I_{k}& -\frac{1}{\sqrt{2}}I_{k}
  \end{array}
\right),
\end{equation}
where $I_{m}$ is the identity matrix of order $n$. Note that $P'$ is the transposition of $P$. Comprehensive consideration of Eqs.(2.2)-(2.4), thus we employ
\begin{equation*}
P'L(G)P=\left(
  \begin{array}{cc}
  L_{A}& 0\\
    0& L_{S}
  \end{array}
\right),
\end{equation*}
where
\begin{eqnarray*}
L_{A}=\left(
\begin{array}{cc}
L_{V_{0}V_{0}}& \sqrt{2}L_{V_{0}V_{1}}\\
\sqrt{2}L_{V_{1}V_{0}}& L_{V_{1}V_{1}}+L_{V_{1}V_{2}}\\
\end{array}
\right),~
L_{S}=L_{V_{1}V_{1}}-L_{V_{1}V_{2}}.
\end{eqnarray*}

\begin{lem}\textup{\cite{Y.L}}
$L_{A},~L_{S}$ thus obtained above can be used in $P_{L(L_{n})}$, namely
\begin{eqnarray*}
P_{L(L_{n})}\big(x\big)=P_{L_{A}}\big(x\big)P_{L_{S}}\big(x\big).
\end{eqnarray*}
\end{lem}
\begin{lem}\textup{\cite{Gut}}
Gutman and mohar show a new calculation formula, namely
\begin{eqnarray*}
Kf(G)=n\sum_{k=2}^{n}\frac{1}{\mu_{k}} ,
\end{eqnarray*}
where $0=\mu_{1}<\mu_{2}\leq\cdots\leq\mu_{n}(n\geq2)$ are the eigenvalues of $L(G)$.
\end{lem}
\begin{lem}\textup{\cite{H.Y.}}
The number of spanning trees of $G$ can also be concluded to the complexity, which is defined by
\begin{eqnarray*}
\mathscr{L}(G)=\frac{1}{n}\prod_{i=1}^{n}\mu_{i}.
\end{eqnarray*}
\end{lem}
\begin{lem}\textup{\cite{Y.X.}}
Assign $C_{n}$ for a cycle with $n$ vertices. The explicit formula of Laplacian spectrum of $C_{n}$ can be constructed, which is expressed as
\begin{eqnarray*}
S_{p}\big(L(C_{n})\big)=\Big\{2-2\cos\frac{(2\pi i)}{n}\Big|1\leq i\leq n\Big\},~
\end{eqnarray*}
and the Kirchhoff index of $C_{n}$ is given by
\begin{eqnarray*}
Kf\big(C_{n}\big)=\frac{n^{3}-n}{12}.
\end{eqnarray*}
\end{lem}

\section{Main results}\label{sct3}
\ \ \ \ \ In this section, we firstly focus on calculating the systematic formula for the Kirchhoff index of $Q_{n}\big(8,4\big)$. By Lemma 2.1, we apply the notion of the eigenvalues of Laplacian matrix to compute the $Kf\big(Q_{n}(8,4)\big)$. Besides, we determine the complexity of $Q_{n}\big(8,4\big)$, which is made up of the products of degree of all vertices. The scenario is similar for $Q'_{n}\big(8,4\big)$.

Given an $n\times n$ matrix $R$, and put deleting the $i_{1}th,~ i_{2}th,~\cdots ,~i_{k}th$ rows and columns of $R$ are expressed as $R\big[\{i_{1},i_{2},\cdots ,~i_{k}\}\big]$. The vertices of $Q_{n}\big(8,4\big)$ are tabulated in Figure 2. Evidently, $\pi=(1,1')(2,2')\cdots (4n,(4n))')$ is an automorphism of $Q_{n}\big(8,4\big)$ and $v_{0}=\emptyset,~v_{1}=\big\{1,2,3,\cdots, 4n\big\},~v_{2}=\big\{1',2',3',\cdots,(4n)'\big\}.$ Maenwile, we assign $L_{A}\big(Q_{n}(8,4)\big)$ and $L_{S}\big(Q_{n}(8,4)\big)$ as $L_{A}$ and $L_{S}$. Then one can get
\begin{eqnarray}
L_{A}=L_{V_{1}V_{1}}+L_{V_{1}V_{2}},~L_{S}=L_{V_{1}V_{1}}-L_{V_{1}V_{2}}.
\end{eqnarray}
By Equation (3.5), we have
\begin{eqnarray*}
 L_{V_1 V_1}&=&
\left(
  \begin{array}{ccccccccccc}
    3 & -1 & & & & & & & & &\\
    -1 & 2 & -1 & & & & & & & &\\
    & -1 & 2 & -1 & & & & & & &\\
    & & -1 & 3 & -1 & & & & & &\\
    & & & -1 & 3 & -1 & & & & &\\
    & & & & -1 & 2 & -1 & & & &\\
    & & & & & & \ddots & & & &\\
    & & & & & & -1 & 3 & -1 & &\\
    & & & & & & & -1 & 2 & -1 &\\
    & & & & & & & & -1 & 2 & -1\\
    & & & & & & & & & -1 & 3\\
  \end{array}
\right)_{(4n)\times(4n)},
\end{eqnarray*}
\begin{eqnarray*}
 L_{V_1 V_2}&=&
\left(
  \begin{array}{ccccccccccc}
    -1 & & & & & & & & & & -1\\
    &  0 & & & & & & & & &\\
    & &  0 & & & & & & & &\\
    & & &  -1 & & & & & & &\\
    & & & &  0 & & & & & &\\
    & & & & &  0 & & & & &\\
    & & & & & & \ddots & & & &\\
    & & & & & & &  -1 & & & \\
    & & & & & & & &  0 & &\\
    & & & & & & & & &  0 &\\
    -1 & & & & & & & & & & -1\\
  \end{array}
\right)_{(4n)\times(4n)}.
\end{eqnarray*}

Hence,
\begin{eqnarray*}
 L_A&=&
\left(
  \begin{array}{ccccccccccc}
    2 & -1 & & & & & & & & & -1\\
    -1 & 2 & -1 & & & & & & & &\\
    & -1 & 2 & -1 & & & & & & &\\
    & & -1 & 2 & -1 & & & & & &\\
    & & & -1 & 2 & -1 & & & & &\\
    & & & & -1 & 2 & -1 & & & &\\
    & & & & & & \ddots & & & &\\
    & & & & & & -1 & 2 & -1 & &\\
    & & & & & & & -1 & 2 & -1 &\\
    & & & & & & & & -1 & 2 & -1\\
    -1 & & & & & & & & & -1 & 2\\
  \end{array}
\right)_{(4n)\times(4n)},
\end{eqnarray*}
and
\begin{eqnarray*}
 L_S&=&
\left(
  \begin{array}{ccccccccccc}
    4 & -1 & & & & & & & & & 1\\
    -1 & 2 & -1 & & & & & & & &\\
    & -1 & 2 & -1 & & & & & & &\\
    & & -1 & 4 & -1 & & & & & &\\
    & & & -1 & 4 & -1 & & & & &\\
    & & & & -1 & 2 & -1 & & & &\\
    & & & & & & \ddots & & & &\\
    & & & & & & -1 & 4 & -1 & &\\
    & & & & & & & -1 & 2 & -1 &\\
    & & & & & & & & -1 & 2 & -1\\
    1 & & & & & & & & & -1 & 4\\
  \end{array}
\right)_{(4n)\times(4n)}.
\end{eqnarray*}

Assuming that
$0=\alpha_{1}<\alpha_{2}\leq\alpha_{3}\leq\cdots\leq\alpha_{4n}$ are the
roots of $P_{L_{A}}\big(x\big)=0$, and
$0<\beta_{1}\leq\beta_{2}\leq\beta_{3}\leq\cdots\leq\beta_{4n}$ are the
roots of $P_{L_{S}}\big(x\big)=0$. By Lemma 2.2, we immediately have
\begin{eqnarray}
Kf\big(Q_{n}(8,4)\big)=8n\Bigg(\sum_{i=2}^{4n}\frac{1}{\alpha_{i}}+\sum_{j=1}^{4n}\frac{1}{\beta_{j}}\Bigg).
\end{eqnarray}
Next, the prime aim of follows is to calculate
$\sum\limits_{i=2}^{4n}\frac{1}{\alpha_{i}}$ and $\sum\limits_{j=1}^{4n}\frac{1}{\beta_{j}}$.\\

Combing with Lemma 2.1 and Lemma 2.4, one can get
\begin{eqnarray}
Kf\big(Q_{n}(8,4)\big)=8n\Bigg(\frac{16n^{2}-1}{12}+\sum_{j=1}^{4n}\frac{1}{\beta_{j}}\Bigg).
\end{eqnarray}

\begin{lem}
Let $\beta_{j}$ be the eigenvalue of $L_{S}$. And assume that $\mu=15+4\sqrt{14}$, $\nu=15-4\sqrt{14}$.
\begin{eqnarray}
\sum_{j=1}^{4n}\frac{1}{\beta_{j}}=\frac{(-1)^{4n-1}b_{4n-1}}{detL_{S}}.
\end{eqnarray}
\end{lem}

\noindent{\bf Proof.}
$$P_{L_{S}}(x)=det(xI-L_{S})=x^{4n}+b_{1}x^{4n-1}+\cdots+b_{4n-1}x+b_{4n},~b_{4n}\neq0.$$

The key to get Eq.(3.8) is to determine $(-1)^{4n-1}b_{4n-1}$ and $detL_{S}$. We demonstrate the $ith$ order principal submatrix $M_{i} ~and~ M'_{i}$ by from the first $i$ rows and columns of the following matrices $L^{1}_{S}$ and $L^{2}_{S}$ respectively. For $i=1,2,...,4n-1$, we have
\begin{eqnarray*}
 L_S^{1}&=&
\left(
  \begin{array}{ccccccccccc}
    4 & -1 & & & & & & & & &\\
    -1 & 2 & -1 & & & & & & & &\\
    & -1 & 2 & -1 & & & & & & &\\
    & & -1 & 4 & -1 & & & & & &\\
    & & & -1 & 4 & -1 & & & & &\\
    & & & & -1 & 2 & -1 & & & &\\
    & & & & & & \ddots & & & &\\
    & & & & & & -1 & 4 & -1 & &\\
    & & & & & & & -1 & 2 & -1 &\\
    & & & & & & & & -1 & 2 & -1\\
    & & & & & & & & & -1 & 4\\
  \end{array}
\right)_{(4n)\times(4n)},
\end{eqnarray*}
and
\begin{eqnarray*}
 L_S^{2}&=&
\left(
  \begin{array}{ccccccccccc}
    2 & -1 & & & & & & & & &\\
    -1 & 2 & -1 & & & & & & & &\\
    & -1 & 4 & -1 & & & & & & &\\
    & & -1 & 4 & -1 & & & & & &\\
    & & & -1 & 2 & -1 & & & & &\\
    & & & & -1 & 2 & -1 & & & &\\
    & & & & & & \ddots & & & &\\
    & & & & & & -1 & 2 & -1 & &\\
    & & & & & & & -1 & 2 & -1 &\\
    & & & & & & & & -1 & 4 & -1\\
    & & & & & & & & & -1 & 4\\
  \end{array}
\right)_{(4n)\times(4n)}.
\end{eqnarray*}

Applying the results above, we can get two Claims.\\

\noindent{\bf Claim 1.} For $1\leq j\leq 4n$,
\begin{eqnarray*}
q_{j}=
\begin{cases}
\Big(\frac{1}{2}+\frac{9\sqrt{14}}{56}\Big)\cdot\mu^{\frac{j}{4}}+\Big(\frac{1}{2}-\frac{9\sqrt{14}}{56}\Big)\cdot\nu^{\frac{j}{4}},& if~j\equiv 0\ (mod\ 4);\\
\Big(2+\frac{31\sqrt{14}}{56}\Big)\cdot\mu^{\frac{j-1}{4}}+\Big(2-\frac{31\sqrt{14}}{56}\Big)\cdot\nu^{\frac{j-1}{4}},& if~j\equiv 1\ (mod\ 4);\\
\Big(\frac{7}{2}+\frac{53\sqrt{14}}{56}\Big)\cdot\mu^{\frac{j-2}{4}}+\Big(\frac{7}{2}-\frac{53\sqrt{14}}{56}\Big)\cdot\nu^{\frac{j-2}{4}},& if~j\equiv 2\ (mod\ 4);\\
\Big(5+\frac{75\sqrt{14}}{56}\Big)\cdot\mu^{\frac{j-3}{4}}+\Big(5-\frac{75\sqrt{14}}{56}\Big)\cdot\nu^{\frac{j-3}{4}},& if~j\equiv 3\ (mod\ 4).
\end{cases}
\end{eqnarray*}

\noindent{\bf Proof of Claim 1.} Put $q_{j}:=detM_{j}$, with $q_{0}=1$, and $q'_{j}:=detM'_{j}$, with $q_{0}=1$. Then by a direct computing, we can get $q_{1}=4,~q_{2}=7,~q_{3}=10,~q_{4}=33,~q_{5}=122,~q_{6}=211,~q_{7}=300,~q_{8}=989.$ For $4\leq j\leq 4n-1$, then
\begin{eqnarray*}
q_{j}=
\begin{cases}
4q_{j-1}-q_{j-2},& if~j\equiv 0 \ (mod\ 4);\\
4q_{j-1}-q_{j-2},& if~j\equiv 1 \ (mod\ 4);\\
2q_{j-1}-q_{j-2},& if~j\equiv 2 \ (mod\ 4);\\
2q_{j-1}-q_{j-2},& if~j\equiv 3 \ (mod\ 4).
\end{cases}
\end{eqnarray*}

For $1\leq j\leq n-1,$ let $a_{j}=m_{4j}; ~0\leq j\leq n-1,~b_{j}=m_{4j+1},~c_{j}=m_{4j+2},~d_{j}=m_{4j+3}.$ Then we can get $a_{1}=33,~b_{0}=4,~c_{0}=7,~d_{0}=10,~b_{1}=122,~c_{1}=211,~d_{1}=300,$ and for $j\geq 2$ , we have
\begin{eqnarray}
\begin{cases}
a_{j}=4d_{j-1}-c_{j-1};\\
b_{j}=4a_{j}-d_{j-1};\\
c_{j}=2b_{j}-a_{j};\\
d_{j}=2c_{j}-b_{j}.
\end{cases}
\end{eqnarray}

From the first three equations in (3.9), one can get $a_{j}=\frac{2}{13}c_{j}+\frac{2}{13}c_{j-1}$. Next, substituting $a_{j}$ to the third equation, one has $b_{j}=\frac{15}{26}c_{j}+\frac{1}{26}c_{j-1}$. Then substituting $b_{j}$ to the fourth equation, we have $d_{j}=\frac{37}{26}c_{j}-\frac{1}{26}c_{j-1}.$ Finally, Substituting $a_{j}$ and $d_{j}$ to the first equation, one has $c_{j}-30c_{j-1}+c_{j=2}=0.$ Thus
\begin{eqnarray*}
c_{j}=k_{1}\cdot\mu^{j}+k_{2}\cdot\nu^{j}.
\end{eqnarray*}

According to $c_{0}=7,c_{1}=211,$ we have
\begin{eqnarray*}
\begin{cases}
k_{1}+k_{2}=7;\\
k_{1}(15+4\sqrt{14})+k_{2}(15-4\sqrt{14})=211.
\end{cases}
\end{eqnarray*}
and
\begin{eqnarray*}
\begin{cases}
k_{1}=(\frac{7}{2}+\frac{53\sqrt{14}}{56}\Big);\\
k_{2}=(\frac{7}{2}-\frac{53\sqrt{14}}{56}\Big).
\end{cases}
\end{eqnarray*}

Thus 
\begin{eqnarray*}
\begin{cases}
a_{j}=\Big(\frac{1}{2}+\frac{9\sqrt{14}}{56}\Big)\cdot\mu^{j}+\Big(\frac{1}{2}-\frac{9\sqrt{14}}{56}\Big)\cdot\nu^{j};\\
b_{j}=\Big(2+\frac{31\sqrt{14}}{56}\Big)\cdot\mu^{j}+\Big(2-\frac{31\sqrt{14}}{56}\Big)\cdot\nu^{j};\\
c_{j}=\Big(\frac{7}{2}+\frac{53\sqrt{14}}{56}\Big)\cdot\mu^{j}+\Big(\frac{7}{2}-\frac{53\sqrt{14}}{56}\Big)\cdot\nu^{j};\\
d_{j}=\Big(5+\frac{75\sqrt{14}}{56}\Big)\cdot\mu^{j}+\Big(5-\frac{75\sqrt{14}}{56}\Big)\cdot\nu^{j},
\end{cases}
\end{eqnarray*}
as desired.\hfill\rule{1ex}{1ex}\\

Using a similar method of Claim 1, we can prove Claim 2.\\

\noindent{\bf Claim 2.} For $1\leq j\leq 4n$,
\begin{eqnarray*}
q^{'}_{j}=
\begin{cases}
\Big(\frac{1}{2}+\frac{11\sqrt{14}}{56}\Big)\cdot\mu^{\frac{j}{4}}+\Big(\frac{1}{2}-\frac{11\sqrt{14}}{56}\Big)\cdot\nu^{\frac{j}{4}},& if~j\equiv 0\ (mod\ 4);\\
\Big(1+\frac{17\sqrt{14}}{56}\Big)\cdot\mu^{\frac{j-1}{4}}+\Big(1-\frac{17\sqrt{14}}{56}\Big)\cdot\nu^{\frac{j-1}{4}},& if~j\equiv 1\ (mod\ 4);\\
\Big(\frac{3}{2}+\frac{23\sqrt{14}}{56}\Big)\cdot\mu^{\frac{j-2}{4}}+\Big(\frac{3}{2}-\frac{23\sqrt{14}}{56}\Big)\cdot\nu^{\frac{j-2}{4}},& if~j\equiv 2\ (mod\ 4);\\
\Big(5+\frac{75\sqrt{14}}{56}\Big)\cdot\mu^{\frac{j-3}{4}}+\Big(5-\frac{75\sqrt{14}}{56}\Big)\cdot\nu^{\frac{j-3}{4}},& if~j\equiv 3\ (mod\ 4).
\end{cases}
\end{eqnarray*}

By the related properties of determinant, it is evidently that
\begin{eqnarray*}
 \det L_S&=&
\left|
  \begin{array}{ccccccc}
    4 & -1 & 0 & \cdots & 0 & 0 & 1\\
    -1 & 2 & -1 & \cdots & 0 & 0 & 0\\
     0 & -1 & 2 & \cdots & 0 & 0 & 0\\
    \vdots & \vdots & \vdots & \ddots & \vdots & \vdots & \vdots\\
    0 & 0 & 0 & \cdots & 2 & -1 & 0\\
    0 & 0 & 0 & \cdots & 0 & 2 & -1\\
    1 & 0 & 0 & \cdots & 0 & -1 & 4  \\
  \end{array}
\right|_{4n}\\
&=&\left|
  \begin{array}{ccccccc}
    4 & -1 & 0 & \cdots & 0 & 0 & 0\\
    -1 & 2 & -1 & \cdots & 0 & 0 & 0\\
     0 & -1 & 2 & \cdots & 0 & 0 & 0\\
    \vdots & \vdots & \vdots & \ddots & \vdots & \vdots & \vdots\\
    0 & 0 & 0 & \cdots & 2 & -1 & 0\\
    0 & 0 & 0 & \cdots & 0 & 2 & -1\\
    1 & 0 & 0 & \cdots & 0 & -1 & 4  \\
  \end{array}
\right|_{4n}+\left|
  \begin{array}{ccccccc}
    4 & -1 & 0 & \cdots & 0 & 0 & 1\\
    -1 & 2 & -1 & \cdots & 0 & 0 & 0\\
     0 & -1 & 2 & \cdots & 0 & 0 & 0\\
    \vdots & \vdots & \vdots & \ddots & \vdots & \vdots & \vdots\\
    0 & 0 & 0 & \cdots & 2 & -1 & 0\\
    0 & 0 & 0 & \cdots & 0 & 2 & 0\\
    1 & 0 & 0 & \cdots & 0 & -1 & 0\\
  \end{array}
\right|_{4n}\\
&=&q_{4n}+(-1)^{4n+1}\cdot(-1)^{4n-1}+(-1)^{4n+1}\cdot\Big[(-1)^{4n-1}+(-1)^{4n}q^{'}_{4n-2}\Big]\\
&=&q_{4n}-q^{'}_{4n-2}+2.\\
\end{eqnarray*}
Together with Claims 1 and 2, we can get one interesting Claim (Claim 3).\\

\noindent{\bf Claim 3.} $detL_{S}=\mu^{n}+\nu^{n}+2.$\\

So far, we have calculated the $detL_{S}$ in Eq.(3.8). Then the prime of the rest is to calculate $(-1)^{4n-1}b_{4n-1}$.\\

\noindent{\bf Claim 4.} $(-1)^{4n-1}b_{4n-1}=\frac{9n\sqrt{14}\big(\mu^{n}-\nu^{n}\big)}{14}.$\\

\noindent{\bf Proof of Claim 4.} We apply the notion that $(-1)^{4n-1}b_{4n-1}$ is the total of all the principal minors of (4$n$-1)\-/order matix $L_{S}$, thus we obtain
\begin{eqnarray*}
(-1)^{4n-1}b_{4n-1}&=&\sum_{j=1}^{4n}detL_{S}[i]\\
&=&\sum_{j=4,j\equiv0(mod\ 4)}^{4n}detL_{S}[j]+\sum_{j=1,j\equiv1(mod\ 4)}^{4n-3}detL_{S}[j]\\
&&+\sum_{j=2,j\equiv2(mod\ 4)}^{4n-2}detL_{S}[j]+\sum_{j=3,j\equiv3(mod\ 4)}^{4n-1}detL_{S}[j].
\end{eqnarray*}

Let
$\left(
\begin{array}{cc}
O & P\\
S & T\\
\end{array}
\right)$
be a block matrix of $L_{S}[j]$
and
N=$\left(
\begin{array}{cc}
0 & -I_{j-1}\\
I_{4n-j} & 0\\
\end{array}
\right)$. It's easy to get
\begin{eqnarray}
N'L_{S}[j]N=\left(
\begin{array}{cc}
0 & -I_{j-1}\\
-I_{4n-j} & 0\\
\end{array}
\right)'
\left(
\begin{array}{cc}
O & P\\
S & T\\
\end{array}
\right)
\left(
\begin{array}{cc}
0 & -I_{j-1}\\
-I_{4n-j} & 0\\
\end{array}
\right)=
\left(
\begin{array}{cc}
T & -S\\
-P & O\\
\end{array}
\right).
\end{eqnarray}

Combining with Eq.(3.10), based on the different value of $j$, we list the following four cases.\\

{\bf Case 1.} For $ j\equiv0(mod\ 4)$~ and ~$4\leq j\leq 4n-4$.
\begin{eqnarray*}
 \ N'L_S[j]N&=&
\left(
  \begin{array}{ccccccc}
    4 & -1 & 0 & \cdots & 0 & 0 & 0\\
    -1 & 2 & -1 & \cdots & 0 & 0 & 0\\
     0 & -1 & 2 & \cdots & 0 & 0 & 0\\
    \vdots & \vdots & \vdots & \ddots & \vdots & \vdots & \vdots\\
    0 & 0 & 0 & \cdots & 4 & -1 & 0\\
    0 & 0 & 0 & \cdots & -1 & 2 & -1\\
    0 & 0 & 0 & \cdots & 0 & -1 & 2  \\
  \end{array}
\right)_{(4n-1)\times(4n-1)}=q_{4n-1},
\end{eqnarray*}
and
\begin{eqnarray*}
\sum_{j=4,j\equiv0(mod\ 4)}^{4n}detL_{S}[j]=nq_{4n-1}=\frac{5n\sqrt{14}\big(\mu^{n}-\nu^{n}\big)}{56}.
\end{eqnarray*}

{\bf Case 2.} For $ j\equiv1(mod\ 4)$~ and ~$1\leq j\leq 4n-3$.
\begin{eqnarray*}
 \ N'L_S[j]N&=&
\left(
  \begin{array}{ccccccc}
    2 & -1 & 0 & \cdots & 0 & 0 & 0\\
    -1 & 2 & -1 & \cdots & 0 & 0 & 0\\
     0 & -1 & 4 & \cdots & 0 & 0 & 0\\
    \vdots & \vdots & \vdots & \ddots & \vdots & \vdots & \vdots\\
    0 & 0 & 0 & \cdots & 2 & -1 & 0\\
    0 & 0 & 0 & \cdots & -1 & 2 & -1\\
    0 & 0 & 0 & \cdots & 0 & -1 & 4  \\
  \end{array}
\right)_{(4n-1)\times(4n-1)}=q'_{4n-1},
\end{eqnarray*}
and
\begin{eqnarray*}
\sum_{j=1,j\equiv1(mod\ 4)}^{4n-3}detL_{S}[j]=nq'_{4n-1}=\frac{5n\sqrt{14}\big(\mu^{n}-\nu^{n}\big)}{56}.
\end{eqnarray*}

{\bf Case 3.} For $ j\equiv2(mod\ 4)$~ and ~$2\leq j\leq 4n-2$.
\begin{eqnarray*}
 \ N'L_S[j]N&=&
\left(
  \begin{array}{ccccccc}
    2 & -1 & 0 & \cdots & 0 & 0 & 0\\
    -1 & 4 & -1 & \cdots & 0 & 0 & 0\\
     0 & -1 & 4 & \cdots & 0 & 0 & 0\\
    \vdots & \vdots & \vdots & \ddots & \vdots & \vdots & \vdots\\
    0 & 0 & 0 & \cdots & 2 & -1 & 0\\
    0 & 0 & 0 & \cdots & -1 & 4 & -1\\
    0 & 0 & 0 & \cdots & 0 & -1 & 4  \\
  \end{array}
\right)_{(4n-1)\times(4n-1)},
\end{eqnarray*}
and
\begin{eqnarray*}
\sum_{j=2,j\equiv2(mod\ 4)}^{4n-2}detL_{S}[j]=n\Big[2\big(4q_{4n-3}-q'_{4n-4}\big)-q_{4n-3}\Big]=\frac{13n\sqrt{14}\big(\mu^{n}-\nu^{n}\big)}{56}).
\end{eqnarray*}

{\bf Case 4.} For $ j\equiv3(mod\ 4)$~ and ~$3\leq j\leq 4n-1$.
\begin{eqnarray*}
 \ N'L_S[j]N&=&
\left(
  \begin{array}{ccccccc}
    4 & -1 & 0 & \cdots & 0 & 0 & 0\\
    -1 & 4 & -1 & \cdots & 0 & 0 & 0\\
     0 & -1 & 2 & \cdots & 0 & 0 & 0\\
    \vdots & \vdots & \vdots & \ddots & \vdots & \vdots & \vdots\\
    0 & 0 & 0 & \cdots & 4 & -1 & 0\\
    0 & 0 & 0 & \cdots & -1 & 4 & -1\\
    0 & 0 & 0 & \cdots & 0 & -1 & 2  \\
  \end{array}
\right)_{(4n-1)\times(4n-1)},
\end{eqnarray*}
and
\begin{eqnarray*}
\sum_{j=3,j\equiv3(mod\ 4)}^{4n-1}detL_{S}[j]=n(4q_{4n-2}-q'_{4n-3})=\frac{13n\sqrt{14}\big(\mu^{n}-\nu^{n}\big)}{56}.
\end{eqnarray*}

Thus, one has the following equation
\begin{eqnarray*}
 (-1)^{4n-1}b_{4n-1}=\frac{9n\sqrt{14}\big(\mu^{n}-\nu^{n}\big)}{14},
\end{eqnarray*}
as desired.\hfill\rule{1ex}{1ex}\\

According to the results of facts 3 and 4, we immediately get
\begin{eqnarray}
\sum_{j=1}^{4n}\frac{1}{\beta_{j}}=\frac{9n\sqrt{14}}{14}\Bigg(\frac{(\mu^{n}-\nu^{n}\big)}{\mu^{n}+\nu^{n}\big)-2}\Bigg).
\end{eqnarray}
\hfill\rule{1ex}{1ex}\\

\begin{thm} Suppose $Q_{n}\big(8,4\big)$ as a M\"{o}bius graph constructed by $n$ octagonals and $n$ quadrilaterals. Then
\begin{eqnarray*}
Kf\big(Q_{n}(8,4)\big)=\frac{32n^{3}-2n}{3}+\frac{36n^{2}\sqrt{14}}{7}\Bigg(\frac{\big(\mu^{n}-\nu^{n}\big)}{\big(\mu^{n}+\nu^{n}\big)+2}\Bigg).
\end{eqnarray*}
\end{thm}

\noindent{\bf Proof.} Substituting Eqs.(3.7) and (3.11) into (3.6), the Kirchhoff index of $Q_{n}(8,4)$ can be expressed
\begin{eqnarray*}
Kf\big(Q_{n}(8,4)\big)&=&8n\Bigg(\sum_{i=2}^{4n}\frac{1}{\alpha_{i}}+\sum_{j=1}^{4n}\frac{1}{\beta_{j}}\Bigg)\\
&=&8n\Bigg(\frac{16n^{2}-1}{12}+\frac{9n\sqrt{14}}{14}\cdot\frac{\big(\mu^{n}-\nu^{n}\big)}{\big(\mu^{n}+\nu\big)+2}\Bigg)\\
&=&\frac{32n^{3}-2n}{3}+\frac{36n^{2}\sqrt{14}}{7}\Bigg(\frac{\big(\mu^{n}-\nu^{n}\big)}{\big(\mu^{n}+\nu\big)+2}\Bigg).
\end{eqnarray*}
The result as desired.\hfill\rule{1ex}{1ex}\\

Next, we aim to compute the Kirchhoff index for $Q'_{n}\big(8,4\big)$. By Figure 2, it is apparently to see $\pi=(1,1')(2,2') \cdots(4n,(4n)')$ is an automorphism of $Q'_{n}\big(8,4\big)$. In other words, $v_{0}=\emptyset,~v_{1}=\big\{1,2,3,\cdots, 4n\big\}~and ~v_{2}=\big\{1',2',3',\cdots,(4n)'\big\}.$ In likewise, we express $L_{A}\big(Q'_{n}(8,4)\big)$ and $L_{S}\big(Q'_{n}(8,4)\big)$ as $L'_{A}$ and $L'_{S}$. Thus one can get $L'_{A}=L_{A}$ and
\begin{eqnarray*}
 L'_S&=&
\left(
  \begin{array}{ccccccccccc}
    4 & -1 & & & & & & & & & -1\\
    -1 & 2 & -1 & & & & & & & &\\
    & -1 & 2 & -1 & & & & & & &\\
    & & -1 & 4 & -1 & & & & & &\\
    & & & -1 & 4 & -1 & & & & &\\
    & & & & -1 & 2 & -1 & & & &\\
    & & & & & & \ddots & & & &\\
    & & & & & & -1 & 4 & -1 & &\\
    & & & & & & & -1 & 2 & -1 &\\
    & & & & & & & & -1 & 2 & -1\\
    -1 & & & & & & & & & -1 & 4\\
  \end{array}
\right)_{(4n)\times(4n)},
\end{eqnarray*}

\begin{thm} Suppose $Q'_{n}\big(8,4\big)$ is  cylinder graph of the octagonal-quadrilateral network, then
\begin{eqnarray*}
Kf\big(Q'_{n}(8,4)\big)=\frac{32n^{3}-2n}{3}+\frac{36n^{2}\sqrt{14}}{7}\Bigg(\frac{\big(\mu^{n}-\nu^{n}\big)}{\big(\mu^{n}+\nu^{n}\big)-2}\Bigg).
\end{eqnarray*}
\end{thm}

\noindent{\bf Proof.} Assuming that
$0=\xi_{1}<\xi_{2}\leq\xi_{3}\leq\cdots\leq\xi_{4n}$ are the
roots of $P_{L_{A}}\big(x\big)=0$, and
$0<\lambda_{1}\leq\lambda_{2}\leq\lambda_{3}\leq\cdots\leq\lambda_{4n}$ are the
roots of $P_{L_{S}}\big(x\big)=0$. By Lemma 2.1, we have $S_{p}\big(Q_{n}(8,4)\big)=\big\{\xi_{1},\xi_{2},\cdots,\xi_{4n},\lambda_{1},\lambda_{2},\cdots,\lambda_{4n}\big\}$. Then we straightforwardly get that
\begin{eqnarray*}
Kf\big(Q'_{n}(8,4)\big)&=&8n\Bigg(\sum_{i=2}^{4n}\frac{1}{\xi_{i}}+\sum_{j=1}^{4n}\frac{1}{\lambda_{j}}\Bigg)\\
&=&2Kf\big(C_{4n}\big)+8n\cdot\frac{ (-1)^{4n-1}b'_{4n-1}}{detL'_{S}}.
\end{eqnarray*}

Using a method similar to the previous paper, it is can directly calculate that $detL'_{S}[j]=detL_{S}[j]$ and $b'_{4n-1}=b_{4n-1}$. Note that

\begin{eqnarray*}
 \det L'_S&=&
\left|
  \begin{array}{ccccccc}
    4 & -1 & 0 & \cdots & 0 & 0 & -1\\
    -1 & 2 & -1 & \cdots & 0 & 0 & 0\\
     0 & -1 & 2 & \cdots & 0 & 0 & 0\\
    \vdots & \vdots & \vdots & \ddots & \vdots & \vdots & \vdots\\
    0 & 0 & 0 & \cdots & 2 & -1 & 0\\
    0 & 0 & 0 & \cdots & 0 & 2 & -1\\
    -1 & 0 & 0 & \cdots & 0 & -1 & 4\\
  \end{array}
\right|_{4n}=q_{4n}-q'_{4n-2}-2.
\end{eqnarray*}

Hence,
\begin{eqnarray*}
Kf\big(Q'_{n}(8,4)\big)&=&8n\Bigg(\sum_{i=2}^{4n}\frac{1}{\xi_{i}}+\sum_{j=1}^{4n}\frac{1}{\lambda_{j}}\Bigg)\\
&=&2Kf\big(C_{4n}\big)+8n\frac{ (-1)^{4n-1}b'_{4n-1}}{detL'_{S}}\\
&=&\frac{32n^{3}-2n}{3}+\frac{36n^{2}\sqrt{14}}{7}\Bigg(\frac{\big(\mu^{n}-\nu^{n}\big)}{\big(\mu^{n}+\nu^{n}\big)-2}\Bigg).
\end{eqnarray*}
This completes the proof.\hfill\rule{1ex}{1ex}\\

The Kirchhoff indices of $Q_{n}\big(8,4\big)$ and $Q'_{n}\big(8,4\big)$ for $n=1,\cdot\cdot\cdot,15$ are tabulated Table 1.\\

\begin{table}[htbp]
\setlength{\abovecaptionskip}{0.15cm} \centering\vspace{.3cm}
\setlength{\tabcolsep}{25pt}
\caption{ The Kirchhoff indices of $Q_{1}(8,4),Q_{2}(8,4),...,Q_{15}(8,4)$ and $Q'_{1}(8,4),Q'_{2}(8,4),...,Q'_{15}(8,4).$}
\begin{tabular}{c|c|c|c}
  \hline
  $G$ & $Kf(G)$ & $G$ & $Kf(G)$ \\
  \hline
  $Q_{1}(8,4)$ & $28.00$ & $Q'_{1}(8,4)$ & $30.57$
  \\
  $Q_{2}(8,4)$ & $160.8$ & $Q'_{2}(8,4)$ & $161.14$
  \\
  $Q_{3}(8,4)$ & $459.17$ & $Q'_{3}(8,4)$ & $459.20$
  \\
  $Q_{4}(8,4)$ & $987.88$ & $Q'_{4}(8,4)$ & $987.89$
  \\
  $Q_{5}(8,4)$ & $1811.07$ & $Q'_{5}(8,4)$ & $1811.07$
  \\
  $Q_{6}(8,4)$ & $2992.74$ & $Q'_{6}(8,4)$ & $2992.74$
  \\
  $Q_{7}(8,4)$ & $4596.90$ & $Q'_{7}(8,4)$ & $4596.90$
  \\
  $Q_{8}(8,4)$ & $6687.54$ & $Q'_{8}(8,4)$ & $6687.54$
  \\
  $Q_{9}(8,4)$ & $9328.67$ & $Q'_{9}(8,4)$ & $9328.67$
  \\
  $Q_{10}(8,4)$ & $12584.28$ & $Q'_{10}(8,4)$ & $12584.28$
  \\
  $Q_{11}(8,4)$ & $16518.38$ & $Q'_{11}(8,4)$ & $16518.38$
  \\
  $Q_{12}$(8,4) & $21194.96$ & $Q'_{12}(8,4)$ & $21194.96$
  \\
  $Q_{13}(8,4)$ & $26678.03$ & $Q'_{13}(8,4)$ & $26678.03$
  \\
  $Q_{14}(8,4)$ & $33031.59$ & $Q'_{14}(8,4)$ & $33031.59$
  \\
  $Q_{15}(8,4)$ & $40319.63$ & $Q'_{15}(8,4)$ & $40319.63$
  \\
  \hline
\end{tabular}
\end{table}

In the sequel, the prime aim is to concentrate on calculate the complexity of $Q_{n}\big(8,4\big)$ and $Q'_{n}\big(8,4\big)$.

\begin{thm} Let $Q_{n}\big(8,4\big)$ and $Q'_{n}\big(8,4\big)$ be M\"{o}bius graph and cylinder graph of the linear octagonal-quadrilateral networks, respectively. Then
\begin{eqnarray*}
\mathscr{L}\big(Q_{n}(8,4)\big)=2n\big(\mu^{n}+\nu^{n}+2\big),
\end{eqnarray*}
and
\begin{eqnarray*}
\mathscr{L}\big(Q'_{n}(8,4)\big)=2n\big(\mu^{n}+\nu^{n}-2\big).
\end{eqnarray*}
\end{thm}

\noindent{\bf Proof.}  With initial conditions as $L_{A}=L'_{A}=L(C_{4n}); ~\tau(C_{4n})=4n.$ Based on Lemma 2.3 and Lemma 2.4, one has the following equation.

\begin{eqnarray*}
\prod_{i=2}^{4n}\alpha_{i}=\prod_{i=2}^{4n}\xi_{i}=4n\tau(C_{4n})=16n^{2}.
\end{eqnarray*}

Note that
\begin{eqnarray*}
\prod_{j=1}^{4n}\beta_{j}=detL_{S}=\big(\mu^{n}+\nu^{n}+2\big),
\end{eqnarray*}
\begin{eqnarray*}
\prod_{j=1}^{4n}\lambda_{j}=detL'_{S}=\big(\mu^{n}+\nu^{n}-2\big).
\end{eqnarray*}

Hence,
\begin{eqnarray*}
\mathscr{L}\big(Q_{n}(8,4)\big)=\frac{1}{8n}\prod_{i=2}^{4n}\alpha_{i}\cdot \prod_{j=1}^{4n}\beta_{j}=2n\big(\mu^{n}+\nu^{n}+2\big),
\end{eqnarray*}
and
\begin{eqnarray*}
\mathscr{L}\big(Q'_{n}(8,4)\big)=\frac{1}{8n}\prod_{i=2}^{4n}\xi_{i}\cdot \prod_{j=1}^{4n}\lambda_{j}=2n\big(\mu^{n}+\nu^{n}-2\big).\\
\end{eqnarray*}
\hfill\rule{1ex}{1ex}\\

The complexity of $Q_{n}(8,4)$ and $Q'_{n}(8,4)$ have been computed, which are tabulated in Table 2.
\begin{table}[htbp]
\setlength{\abovecaptionskip}{0.15cm}
  \centering \vspace{.2cm}
  \setlength{\tabcolsep}{10pt}
\caption{The complexity of $Q_{1},Q_{2},...,Q_{8}$ and $Q'_{1},Q'_{2},...,Q'_{8}$.}
\begin{tabular}{c|c|c|c}
  \hline
  ~~~~~~$G$ ~~~~~~&~~~~~~$\mathscr{L}(G)$ ~~~~~~&~~~~~~$G$ ~~~~~~&~~~~~~$\mathscr{L}(G)$~~~~~~ \\
  \hline
  $Q_{1}$ & $64$ & $Q'_{1}$ & $56$ \\
  $Q_{2}$ & $3600$ &  $Q'_{2}$ & $3584$ \\
  $Q_{3}$ & $161472$ & $Q'_{3}$ & $161448$ \\
  $Q_{4}$ & $6451232$ & $Q'_{4}$ & $6451200$ \\
  $Q_{5}$ & $241651520$ & $Q'_{5}$ & $241651480$\\
  $Q_{6}$ & $8689777200$ & $Q'_{6}$ & $8689777152$\\
  $Q_{7}$ & $303803889088$ & $Q'_{7}$ & $303803889032$\\
  $Q_{8}$ & $10404546969664$ & $Q'_{8}$ & $10404546969600$\\
  \hline
\end{tabular}
\end{table}
\section{Relation between Kirchhoff index and Wiener index of $Q_{n}(8,4)$ and $Q'_{n}(8,4)$}\label{sct4}
First part of this section is a certification of the Wiener index of $Q_{n}\big(8,4\big)$ and $Q'_{n}\big(8,4\big)$, and the rest is Wiener index and  Kirchhoff index for comparison.\\
\begin{thm}For the network $Q_{n}\big(8,4\big)$,
\begin{eqnarray*}
W\big(Q_{n}(8,4)\big)=32n^3+16n^2+4n.
\end{eqnarray*}
\end{thm}

\noindent{\bf Proof.}
For the M\"{o}bius graph $Q_{n}\big(8,4\big)$, two types of vertices can be considered $(n\geq2)$.\\

(a) vertex with degree two,\\

(b) vertex with degree three.\\

Hence, for type a, we obtain
\begin{eqnarray}
\omega_{1}(k)=4n\Bigg[\Big(\sum_{k=1}^{2n}k+\sum_{k=1}^{2n}k+\sum_{k=2}^{2n}k+\sum_{k=3}^{2n}k\Big)+3+4\Bigg],
\end{eqnarray}
and for type b, we have
\begin{eqnarray}
\omega_{2}(k)=4n\Bigg(\sum_{k=1}^{2n}k+\sum_{k=1}^{2n}k+\sum_{k=1}^{2n}k+\sum_{k=2}^{2n}k\Bigg).
\end{eqnarray}

Summing the Eqs.(4.12), (4.13) together, we get the Eq.(4.14).\\
\begin{equation}
\begin{split}
W\big(Q_{n}(8,4)\big)&=\frac{\omega_{1}(k)+\omega_{2}(k)}{2}\\
&=\frac{4n\Bigg(\sum\limits_{k=1}^{2n}k+\sum\limits_{k=1}^{2n}k+\sum\limits_{k=1}^{2n}k+\sum\limits_{k=2}^{2n}k+\sum\limits_{k=1}^{2n}k+\sum\limits_{k=1}^{2n}k+\sum\limits_{k=2}^{2n}k+\sum\limits_{k=3}^{2n}k+7\Bigg)}{2}\\
&=32n^3+16n^2+4n.
\end{split}
\end{equation}
\hfill\rule{1ex}{1ex}\\
\begin{thm}
Turning to the cylinder graph $Q'_{n}(8,4)$, then
\begin{eqnarray*}
W\big(Q'_{n}(8,4)\big)=32n^3+16n^2-8n.
\end{eqnarray*}
\end{thm}
\noindent{\bf Proof.}
According to the proof of $W\big(Q_{n}(8,4)\big)$, we have the Eqs.(4.15), (4.16) by similar approach.
\begin{eqnarray}
\omega_{1}(k)=4n\Big[\Big(\sum_{k=1}^{2n-1}k+\sum_{k=1}^{2n}k+\sum_{k=2}^{2n}k+\sum_{k=3}^{2n+1}k\Big)+3+4\Big].
\end{eqnarray}

\begin{eqnarray}
\omega_{2}(k)=4n\Bigg[1+\Big(\sum_{k=1}^{2n-1}k+\sum_{k=1}^{2n}k+\sum_{k=2}^{2n}k+\sum_{k=2}^{2n+1}k\Big)\Bigg].
\end{eqnarray}

Combining the Eqs.(4.15), (4.16) together, we derive the Eq.(4.17).\\
\begin{equation}
\begin{split}
W\big(Q'_{n}(8,4)\big)&=\frac{\omega_{1}(k)+\omega_{2}(k)}{2}\\
&=\frac{4n\Bigg(\sum\limits_{k=1}^{2n-1}k+\sum\limits_{k=1}^{2n}k+\sum\limits_{k=2}^{2n}k+\sum\limits_{k=2}^{2n+1}k+\sum\limits_{k=1}^{2n-1}k+\sum\limits_{k=1}^{2n}k+\sum\limits_{k=2}^{2n}k+\sum\limits_{k=3}^{2n+1}k+8\Bigg)}{2}\\
&=32n^3+16n^2-8n.
\end{split}
\end{equation}
\hfill\rule{1ex}{1ex}\\
Whatever $n$ is even or odd, we deduced the Wiener indices of the M\"{o}bius and cylinder graph of the linear octagonal-quadrilateral networks must be tantamount.\\

\noindent{\bf Corollary 4.3.}
For different values of $n$, we illustrate the ratios of the Wiener index to the Kirchhoff index for both $Q_{n}(8,4)$ and $Q'_{n}(8,4)$ are gradually closing to $3$.\\
\begin{equation*}
\lim_{n\to\infty}\frac{W\big(Q_{n}(8,4)\big)}{Kf\big(Q_{n}(8,4)\big)}=3,
\end{equation*}
and
\begin{equation*}
\lim_{n\to\infty}\frac{W\big(Q'_{n}(8,4)\big)}{Kf\big(Q'_{n}(8,4)\big)}=3.
\end{equation*}

%
\section{Conclusion}\label{sct5}
\ \ \ \ \
A central result of this paper is a systematic expression for the Kirchhoff indices and complexity of the linear M\"{o}bius and cylinder octagonal-quadrilateral networks in terms of the reciprocal of the eigenvalues of the Laplacian matrix. By comparatively analyzing the correlation ability of topological indices with Wiener indices and Kirchhoff indices, we have discovered the underlying relationship between them that the ratio value consistently being three. Finally, we expect the results developed here could become useful part of the research in other polygon chains and their variants.

\section*{Funding}
\ \ \ \ \ This work was supported in part by Anhui Provincial Natural Science Foundation under Grant 2008085J01, and by Natural Science Fund of Education Department of Anhui Province under Grant KJ2020A0478 and National Natural Science Foundation of China Grant 11601006, 12001008.



\end{document}